\documentclass[11pt,a4paper]{article}
\usepackage{enumerate,amsthm,amsmath,graphicx,calc}
\usepackage[lmargin=28mm,rmargin=28mm,tmargin=29mm,bmargin=29mm,footskip=10mm]{geometry}
\usepackage[sort&compress,numbers]{natbib}

\theoremstyle{plain}
\newtheorem{theorem}{Theorem}
\newtheorem{lemma}[theorem]{Lemma}
\newtheorem{corollary}[theorem]{Corollary}

\renewcommand{\baselinestretch}{1.175}
\newcommand{\CC}{\ensuremath{D}}

\newcommand{\fifth}{\ensuremath{\tfrac{1}{5}}}
\newcommand{\half}{\ensuremath{\tfrac{1}{2}}}
\newcommand{\Oh}[1]{\ensuremath{\protect\mathcal{O}(#1)}}
\newcommand{\ceil}[1]{\ensuremath{\protect\lceil#1\rceil}}

\newcommand{\floor}[1]{\ensuremath{\protect\lfloor#1\rfloor}}

\begin{document}

\renewcommand{\thefootnote}{\fnsymbol{footnote}}

\begin{center}
\vspace*{2ex}
{\LARGE\bf Compatible Geometric Matchings\footnote{This work was initiated at the \emph{3rd U.P.C.\ Workshop on Combinatorial Geometry} (Caldes de Malavella, Catalunya, Spain, May 8--12, 2006).}
\\[3ex]}

\renewcommand{\thefootnote}{\arabic{footnote}}

{\large 
Oswin Aichholzer\footnote[1]{Institute for Software Technology, Graz University of Technology, Austria (\texttt{oaich@ist.TUGraz.at}). Supported by the Austrian FWF Joint Research Project 'Industrial Geometry' S9205-N12.}
\hspace*{1em}
Sergey Bereg\footnote[2]{Department of Computer Science, University of Texas at Dallas, U.S.A. (\texttt{besp@utdallas.edu}).}
\hspace*{1em}
Adrian Dumitrescu\footnote[3]{Department of Computer Science, University of Wisconsin-Milwaukee, U.S.A. (\texttt{ad@cs.uwm.edu}). Research partially supported by NSF CAREER grant CCF-0444188.} 
\hspace*{1em}
Alfredo Garc\'ia\footnote[4]{Departamento de M\'etodos Estad\'isticos, Universidad de Zaragoza, Spain (\texttt{olaverri@unizar.es}). Research supported by the project MEC MTM2006-01267.}
\\[1ex]
Clemens Huemer\footnote[5]{Departament de Matem{\`a}tica Aplicada II, Universitat Polit{\`e}cnica de Catalunya, Spain (\texttt{\{ferran.hurtado,clemens.huemer,david.wood\}@upc.edu}). Research supported by the projects MEC MTM2006-01267 and DURSI 2005SGR00692. The research of David Wood is supported by a Marie Curie Fellowship of the European Commission under contract MEIF-CT-2006-023865.}
\hspace*{1em}
Ferran Hurtado\footnotemark[5] 
\hspace*{1em}
Mikio Kano\footnote[6]{Department of Computer and Information Sciences, Ibaraki University, Japan (\texttt{kano@mx.ibaraki.ac.jp}).}
\\[1ex]
Alberto M\'{a}rquez\footnote[7]{Departamento de Matem\'atica Aplicada I, Universidad de Sevilla, Spain (\texttt{almar@us.es}).}
\hspace*{1em}
David Rappaport\footnote[8]{School of Computing, Queen's University, Canada (\texttt{daver@cs.queensu.ca}).
Research supported by NSERC of Canada Discovery Grant 9204.}
\hspace*{1em}
Shakhar Smorodinsky\footnote[9]{Department of Mathematics, Ben-Gurion University,
Israel (\texttt{shakhar@math.bgu.ac.il}).}
\\[1ex]
Diane Souvaine\footnote[10]{Department of Computer Science, Tufts University, U.S.A. (\texttt{dls@cs.tufts.edu}).}
\hspace*{1em}
Jorge Urrutia\footnote[11]{Instituto de Matem\'{a}ticas, Universidad Nacional Aut\'onoma de M\'exico, M\'exico (\texttt{urrutia@math.unam.mx}). Supported by
CONACYT of Mexico, Proyecto SEP-2004-Co1-45876.}
\hspace*{1em}
David R.\ Wood\footnotemark[5]
\\[4ex]
}


\today\\[4ex]

\begin{minipage}{12cm}
\small
\textbf{Abstract:}
This paper studies non-crossing geometric perfect matchings. Two such perfect matchings are \emph{compatible} if they have the same vertex set and their union is also non-crossing. Our first result states that for any two perfect matchings $M$ and $M'$ of the same set of $n$ points, for some $k\in\Oh{\log n}$, there is a sequence of perfect matchings $M=M_0,M_1,\dots,M_k=M'$, such that each $M_i$ is compatible with $M_{i+1}$. This improves the previous best bound of $k\leq n-2$. We then study the conjecture: \emph{every perfect matching with an even number of edges has an edge-disjoint compatible perfect matching}. We introduce a sequence of stronger conjectures that imply this conjecture, and prove the strongest of these conjectures in the case of perfect matchings that consist of vertical and horizontal segments. Finally, we prove that every perfect matching with $n$ edges has an edge-disjoint compatible matching with approximately $4n/5$ edges.
\end{minipage}

\vspace*{2ex}

\end{center}







\newpage
\section{Introduction}


A \emph{geometric graph} is a simple graph $G$, where the vertex-set $V(G)$ is a finite set of points in the plane, and each edge in $E(G)$ is a closed segment whose endpoints belong to $V(G)$. Throughout this paper, we assume that $V(G)$ is in general position; that is, no three vertices are collinear. A geometric graph is \emph{non-crossing} if no two edges cross. That is, two edges may intersect only at a common endpoint. Two non-crossing geometric graphs are \emph{compatible} if they have the same vertex set and their union is non-crossing. 


In this paper, a \emph{matching} is a non-crossing geometric graph in which every vertex has degree at most one. A matching is \emph{perfect} if every vertex has degree exactly one. We say that a (perfect) matching is a (\emph{perfect})  \emph{matching of} its vertex set. Our focus is on compatible perfect matchings. 


We first consider the problem of transforming a given perfect matching into another given perfect matching on the same vertex set. Let $S$ be a set of $n$ points in the plane, with $n$ even. For perfect matchings $M$ and $M'$ of $S$, a \emph{transformation between $M$ and $M'$ of length} $k$ is a sequence $M=M_0,M_1,\dots,M_k=M'$ of perfect matchings of $S$, such that $M_i$ is compatible with $M_{i+1}$, for all $i\in\{0,1,\dots,k-1\}$. \citet{HHNR-GC05} proved that there is a transformation of length $n-2$ between any given pair of perfect matchings of $S$. The first contribution of this paper is to improve this bound from $n-2$ to \Oh{\log n}. This result is proved in Section~\ref{sec:Transforming}.

The remainder of the paper is concerned with the following conjecture. Two geometric graphs are \emph{disjoint} if they have no edge in common. A matching is \emph{even} or \emph{odd} if the number of edges is even or odd. 

\medskip\noindent\textbf{Compatible Matching Conjecture.} 
For every even perfect matching $M$, there is a perfect matching that is disjoint and compatible with $M$.
\medskip

Note that the assumption that the given perfect matching is even is necessary, since there are odd perfect matchings that have no disjoint compatible perfect matching, as described in Section~\ref{sec:Odd}. 

Section~\ref{sec:CompatibleDisjoint} describes progress toward the proof of this conjecture. In particular, we introduce a sequence of stronger conjectures that imply the Compatible Matching Conjecture. 

In the next two sections we establish the Compatible Matching Conjecture for the following special cases: perfect matchings that consist of vertical and horizontal segments (Section~\ref{sec:VerticalHorizontal}), and  perfect matchings that arise from convex-hull-connected sets of segments (Section~\ref{CHC}).

In the final two sections we consider two relaxations of the Compatible Matching Conjecture. First we relax the requirement that the matching is perfect, and we prove that every perfect matching with $n$ edges has a disjoint compatible (partial) matching with approximately $4n/5$ edges (Section~\ref{ThreeQuarters}). Finally, we prove a weakened version of the Compatible Matching Conjecture in which certain types of crossings are allowed (Section~\ref{WithCrossings}).

\subsection{Related Work}

Instead of transforming perfect matchings, \citet{AAH-CGTA02} considered transforming spanning trees of a fixed set of $n$ points, and established the following results. Start with any non-crossing spanning tree $T$, and let $f(T)$ be the shortest spanning tree that does not cross $T$. Then $f(T)$ is non-crossing. In addition, iterating the operator $f$ must stop at some point, because the total length of the edges is decreasing. \citet{AAH-CGTA02} proved that this process always leads to a minimum spanning tree, for every starting tree $T$. Moreover, it takes \Oh{\log n} steps to reach a minimum spanning tree, and for some starting trees, $\Omega(\log n)$ steps are required. As a corollary, there is a transformation of length \Oh{\log n} between any two spanning trees. Whether this bound is tight is of some interest. Partially motivated by connections with pseudo-triangulations, \citet{AAHK-IPL06} conjectured that there is a transformation of length $o(\log n)$ between any two spanning trees. Recently \citet{BRUW} proved an $\Omega(\log n/\log\log n)$ lower bound for this question.

There is another problem that has attracted substantial research and is closely related to the Compatible Matching Conjecture. In general, given a set $S$ of $n$ pairwise disjoint segments it is not always possible to form a polygon with $2n$ sides such that every second segment on its boundary belongs to $S$ (an \emph{alternating} polygon). Toussaint raised the computational problem of deciding whether an alternating polygon exists, which was extensively studied by Rappaport and other authors \citep{Rappaport-SJC89,RIT-DCG90}. Later \citet{M-CGTA92} conjectured that there is a polygon such that every segment from $S$ is a side or an internal diagonal (a \emph{circumscribing polygon}); this was disproved by \citet{UW-CGTA92}. \citet{PR-CGTA98} proved that there is a circumscribing polygon of size $\Omega(n^{1/3})$ (although this cycle may cross the other segments). \citet{M-CGTA92} also conjectured that there is a polygon such that every segment from $S$ is a side, an internal diagonal, or an external diagonal. This conjecture was finally proved by \citet{HT-CGTA03}.

The Compatible Matching Conjecture follows the lines of the original formulation of the preceding problem, as it implies that there is a \emph{set} of pairwise disjoint simple polygons, with a total of $2n$ edges, such that every segment from $S$ lies on the boundary of one of them. 


\section{Tools}

\subsection{Matchings in a Polygon} 

The following result by \citet{AGHTU05} is used repeatedly throughout the paper.

\begin{lemma}[\citep{AGHTU05}]
\label{lem:reflex}
Let $P$ be a simple polygon, let $R$ be the set of reflex vertices of $P$, and let $S$ be any finite set of points on the boundary of $P$ or in its interior, such that $R\subseteq S$ and $|S|$ is even. Then there is a perfect matching $M$ of $S$ such that every segment in $M$ is inside the (closed) polygon $P$.
\end{lemma}

While in general, the Compatible Matching Conjecture is false for odd perfect matchings of point sets in convex position, the following lemma provides an important special case when a disjoint compatible perfect matching always exists. 

\begin{lemma}
\label{lem:Convex}
Let $P$ be a set of points in convex position. Let $M$ be a matching of $P$ such that every segment in $M$ is on the boundary of the convex hull of $P$. Then there is a perfect matching of $P$ that is disjoint and compatible with $M$ if and only if $|P|$ is even and if $|P|=2$ then $E(M)=\emptyset$.
\end{lemma}

\begin{proof}
The necessity of the conditions are obvious. We prove the sufficiency by induction on $|P|$. The base cases with $|P|\leq4$ are easily verified. Now suppose that $|P|\geq6$ is even. Thus there are consecutive vertices $v$ and $w$ in $P$ that are not adjacent in $M$. Let $P':=P-\{v,w\}$. Let $M'$ be the subgraph of $M$ induced by $P'$. Thus $M'$ is a matching of $P'$ such that every segment in $M'$ is on the convex hull of $P'$. Since $|P'|\geq4$ is even, by induction, $M'$ has a disjoint compatible perfect matching $M''$. Let $M'''$ be the geometric graph obtained from $M''$ by adding the vertices $v$ and $w$, and adding the edge $vw$. Now $M'''$ is non-crossing, since $v$ and $w$ are consecutive on the convex hull of $P$. Since $vw\not\in E(M)$, $M$ and $M'''$ are disjoint.
\end{proof}



\subsection{Segment Extensions} 
\label{sec:Extensions}


Let $M$ be a perfect matching, and let $C$ be a (possibly unbounded) convex set in the plane, such that every segment in $M$ that intersects $C$ has at least one endpoint in $C$. 
Let $M_1$ be the set of segments in $M$ with exactly one endpoint in $C$. Let $M_2$ be the set of segments in $M$ with both endpoints in $C$. We ignore the segments in $M$ outside of $C$. 

\begin{figure}[htb]
\centering\includegraphics[scale=0.5]{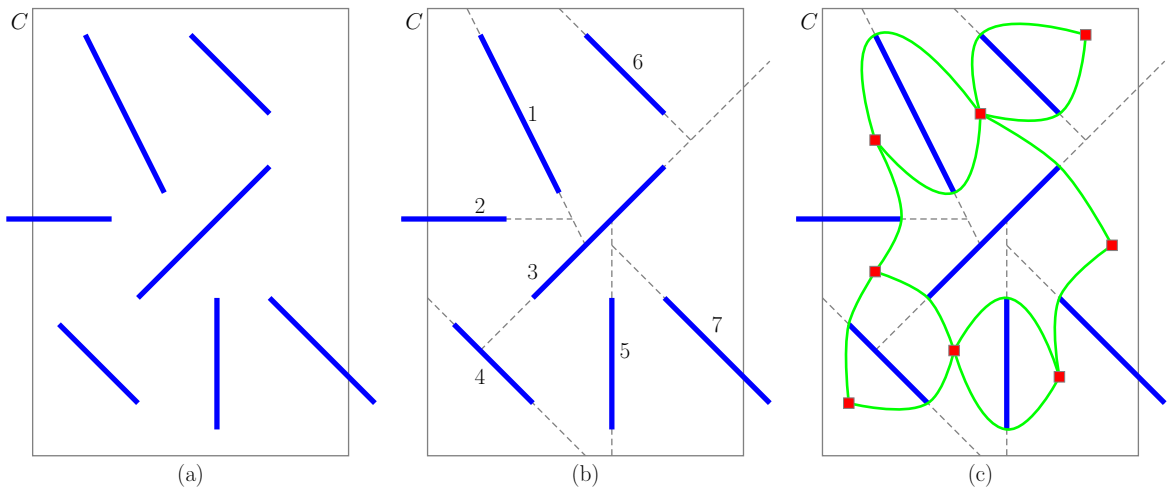}
\caption{(a) A given perfect matching $M$ and convex region $C$. (b) The extension of $M$ in the order shown produces a subdivision of $C$. (c) The associated dual multigraph $G$.}
\label{fig:Extension} 
\end{figure}

As illustrated in Figure~\ref{fig:Extension}, an \emph{extension} of $M$ in $C$ is a set of segments and rays obtained as follows. For each segment $s\in M_1\cup M_2$ in some given order, extend $s$ by a ray, in both directions if $s\in M_2$, and in the direction into $C$ if $s\in M_1$. Each ray is extended until it hits another segment, the boundary of $C$, or a previous extension, or the ray \emph{goes to infinity} if it is not blocked. An extension $L$ of $M$ defines a convex subdivision of $C$ with $|M_1|+|M_2|+1$ cells, since the extension of each segment splits one cell into two cells. 

Since the vertices of $M$ are in general position by assumption, each vertex of $M$ that is in $C$ is on the boundary of exactly two cells of the convex subdivision. The \emph{dual multigraph} $G$ of $L$ is the (non-geometric) multigraph whose vertices are the cells of this convex subdivision. For every vertex $v$ of $M$ that is in $C$, add an edge to $G$ between the vertices that correspond to the two cells of the convex subdivision of which $v$ is on the boundary. Thus $G$ has $|M_1|+|M_2|+1$ vertices and $|M_1|+2|M_2|$ edges. Since $G$ is obtained by a series of vertex splitting\footnote{Let $v$ be a vertex in a connected graph $G$. Let $S$ be a subset of the neighbours of $v$. Let $G'$ be the graph obtained from $G$ by deleting the edges from $v$ to $S$, and introducing a new vertex $v'$ adjacent to $v$ and to each vertex in $S$. Then $G'$ is said to be obtained from $G$ by \emph{splitting} $v$. Clearly $G'$ is also connected.} operations, $G$ is connected. 

The above properties of extensions of perfect matchings are folklore \citep{BHT-DCG01,EHKN-CGTA00,ORourke87}.

\subsection{Even Orientations} 

Our third tool is non-geometric. A \emph{multigraph} allows parallel edges. A multigraph is \emph{even} if it has an even number of edges, and \emph{odd} otherwise. An orientation of a multigraph is \emph{even} if every vertex has even indegree. \citet{FJS-DAM01} and \citet{FrankKiraly-Comb02} characterised when a multigraph admits an even orientation.

\begin{lemma}[\citep{FJS-DAM01,FrankKiraly-Comb02}]
\label{lem:EvenOrientation}
A multigraph admits an even orientation if and only if every component is even.
\end{lemma}

Note that a tree $T$ with an even number of edges has a unique even orientation. In particular, for each edge $vw$ of $T$, consider the subtrees $T_v$ and $T_w$ obtained from $T$ by deleting $vw$, where $v\in V(T_v)$ and $w\in V(T_w)$. Then one of $T_v$  and $T_w$ has an even number of edges and the other has an odd number of edges. Say $|E(T_v)|$ is even. Consider an even orientation of $T$. Then $vw$ is oriented from $v$ to $w$, as otherwise $T_v$ plus the edge $vw$ would be evenly oriented, but this subtree has an odd number of edges, which is clearly impossible. Similarly, if $|E(T_w)|$ is even then $vw$ is oriented from $w$ to $v$ in every even orientation. Conversely, if we orient each edge $vw$ as described above, then it is easily seen that we obtain an even orientation of $T$.    

\section{Transforming Matchings}
\label{sec:Transforming}

In this section we prove the following theorem.

\begin{theorem}
\label{thm:Transform}
For every set $S$ of $2n$ points in general position, there is a transformation of length at most $2\ceil{\log_2n}$ between any given pair of perfect matchings of $S$. 
\end{theorem}


We begin with some preliminary lemmas.

\begin{lemma}
\label{lem:HalfPlane}
Let $M$ be a perfect matching. Let $t$ be a line cutting an even number of segments in $M$, but containing no vertex of $M$. Let $H$ be a halfplane determined by $t$. Let $S$ be the set of vertices of $M$ that are in $H$. Then there is a perfect matching $M'$ of $S$ such that $M\cup M'$ is non-crossing.
\end{lemma}

\begin{proof}[First Proof] 
Say $m$ segments of $M$ are cut by $t$, and $n$ segments of $M$ are contained in $H$. As described in Section~\ref{sec:Extensions}, consider an extension of $M$ in $H$. The obtained subdivision of $H$ has $m+n+1$ convex cells, and the dual multigraph $G$ is connected. Since $m$ is even, the number of edges of $G$, $m+2n$, is also even. By Lemma~\ref{lem:EvenOrientation}, $G$ admits an even orientation. Thus each vertex in $S$ can be assigned to one of its two adjacent cells, so that each cell $C$ is assigned an even number of vertices. Let $S_C$ be the set of vertices assigned to cell $C$. Since $C$ is convex, there is a perfect matching of $S_C$ that is compatible with the matching of $S_C$ induced by $M$. (We cannot conclude that these matchings are disjoint, as in Lemma~\ref{lem:Convex}, since it is possible that $|S_C|=2$ and the two points are endpoints of the same segment.)\ The union of these matchings, taken over all the convex cells $C$, is a perfect matching $M'$ of $S$, such that $M\cup M'$ is non-crossing.
\end{proof}

\begin{figure}[htb]
\centering\includegraphics{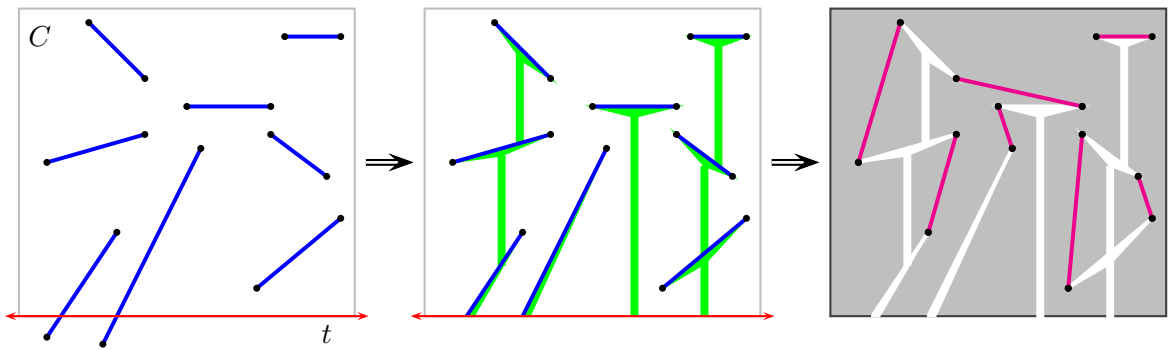}
\caption{Second proof of Lemma~\ref{lem:HalfPlane}.}
\label{fig:EvenCut} 
\end{figure}

\begin{proof}[Second Proof]
Without loss of generality, $t$ is horizontal, and no segment in $M$ is vertical. Let $C$ be a rectangle containing $S$ whose base side is contained in $t$. For each segment $vw$ of $M$ with at least one endpoint in $S$, let $x$ be a point infinitesimally below the midpoint of $vw$. Now, thicken $vw$ into the triangle $vxw$. Moreover, if both $v$ and $w$ are in $S$, then draw an infinitesimally wide axis-parallel rectangle from $x$ downward until it reaches $t$ or another segment of $M$. As illustrated in Figure~\ref{fig:EvenCut}, removing the infinitesimal elements from $C$, we obtain a simple polygon $P$ whose reflex vertices are precisely the vertices in $S$. By Lemma~\ref{lem:reflex} with $R=S$, there is a perfect matching $M'$ of $S$, such that every segment in $M'$ is inside $P$. Thus $M\cup M'$ is non-crossing.
\end{proof}

\begin{lemma}
\label{lem:EvenCut}
Let $M$ be a perfect matching. Let $t$ be a line cutting an even number of segments of $M$, but containing no vertex of $M$. Let $S_1$ and $S_2$ be the sets of vertices of $M$ lying in the two halfplanes determined by $t$. Then there is a perfect matching $M_1$ of $S_1$ and a perfect matching $M_2$ of $S_2$, such that $M$ and $M_1\cup M_2$ are compatible (but not necessarily disjoint).
\end{lemma}

\begin{proof}
Apply Lemma~\ref{lem:HalfPlane} to each half-plane determined by $t$, to obtain a perfect matching $M_1$ of $S_1$ and a perfect matching $M_2$ of $S_2$, such that $M\cup M_1$ is non-crossing and $M\cup M_2$ is non-crossing. Now $M_1\cup M_2$ is non-crossing since $M_1$ and $M_2$ are separated by $t$. Hence $M_1\cup M_2$ is a perfect matching of the vertex set of $M$, and it is compatible with $M$.
\end{proof}

Let $S$ be a set of $2n$ points in general position in the plane. Without loss of generality, no two points in $S$ have the same X-coordinate. Order the points $p_1,p_2,\dots,p_{2n}$ in increasing order of their X-coordinates. Let $N(S)$ be the \emph{canonical} perfect matching that links $p_{2i-1}$ with $p_{2i}$ for each $i\in\{1,2,\dots,n\}$. 


\begin{lemma}
\label{lem:MakeCanonical}
For every set $S$ of $2n$ points in general position, and for every perfect matching $M$ of $S$, there is a transformation between $M$ and $N(S)$ of length $\ceil{\log_2n}$.
\end{lemma}


\begin{proof}
Let $k(n)=\ceil{\log_2n}$. We proceed by induction on $n$. With $n=1$, every perfect matching of $S$ is canonical, and we are done since $k(1)=0$. Now assume that $n>1$ and the lemma holds for all values less than $n$. Let $t$ be a vertical line with $2\floor{n/2}$ points of $S$ to the left of $t$, and $2\ceil{n/2}$ points of $S$ to the right of $t$. Let $S^{\ell}$ and $S^r$ be the sets of points in $S$ respectively to the left and right of $t$. Say $t$ cuts $m$ edges of $M$. The $2\floor{n/2}-m$ points of $S^{\ell}$ that are incident to an edge of $M$ not cut by $t$ are matched by $M$. Thus $m$ is even. By Lemma~\ref{lem:EvenCut}, there is a perfect matching $M^{\ell}$ of $S^{\ell}$ and a perfect matching $M^r$ of $S^r$, such that $M$ and $M^{\ell}\cup M^r$ are compatible. Now apply induction to $M^\ell$ and $M^r$. Observe that $k(\floor{n/2})\leq k(\ceil{n/2})\leq k(n)-1$. Thus there are transformations $$M^\ell=M^\ell_0,M^\ell_1,\dots,M^\ell_{k(n)-1}=N(S^\ell)\text{ and }M^r=M^r_0,M^r_1,\dots,M^r_{k(n)-1}=N(S^r).$$
Hence each $M^\ell_i$ is compatible with $M^\ell_{i+1}$, and each $M^r_i$ is compatible with $M^r_{i+1}$. Let $M_i:=M^\ell_{i-1}\cup M^r_{i-1}$ for each $i\in\{1,2,\dots,k(n)\}$. Since 
$M^\ell_{i-1}$ and $M^r_{i-1}$ are separated by $t$, $M_i$ is a perfect matching of $S$, and $M_i$ is compatible with $M_{i+1}$ for each $i\in\{1,2,\dots,k(n)-1\}$. By Lemma~\ref{lem:EvenCut}, $M$ and $M_1=M^{\ell}\cup M^r$ are compatible. Since $N(S)=N(S^\ell)\cup N(S^r)=M_{k(n)}$, $$M,M_1,\dots,M_{k(n)}$$ is a transformation between $M$ and $N(S)$ of length $k(n)$.
\end{proof}



\begin{proof}[Proof of Theorem~\ref{thm:Transform}] 
For perfect matchings $M$ and $M'$ of $S$, by Lemma~\ref{lem:MakeCanonical}, there are transformations $$M=M_0,M_1,\dots,M_{k(n)}=N(S)\text{  and }M'=M'_0,M'_1,\dots,M'_{k(n)}=N(S).$$ Thus $M=M_0,M_1,\dots,M_{k(n)},M'_{k(n)-1},M'_{k(n)-2}\dots,M'_0=M'$ is a transformation between $M$ and $M'$ of length $2k(n)$.
\end{proof}

\section{Odd Matchings}
\label{sec:Odd}

\begin{figure}[htb]
\centering\includegraphics{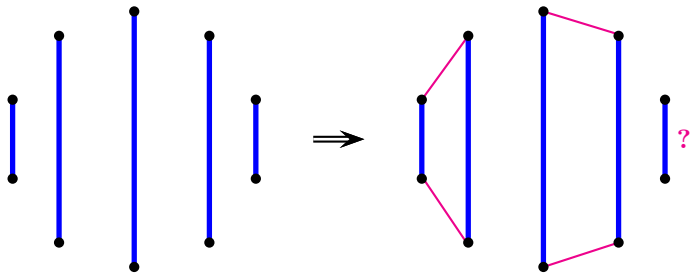}
\caption{A 5-edge perfect matching with no disjoint compatible perfect matching.}
\label{OddExample}
\end{figure}

In the remainder of the paper we study the Disjoint Compatible Matching conjecture. First, in this section, we show why this conjecture is false for odd perfect matchings. That is, we describe classes of odd perfect matchings that have no disjoint compatible perfect matching. It is easily seen that an odd number of parallel chords of a circle form such a matching, as illustrated in Figure~\ref{OddExample}. 

\begin{figure}[htb]
\centering\includegraphics{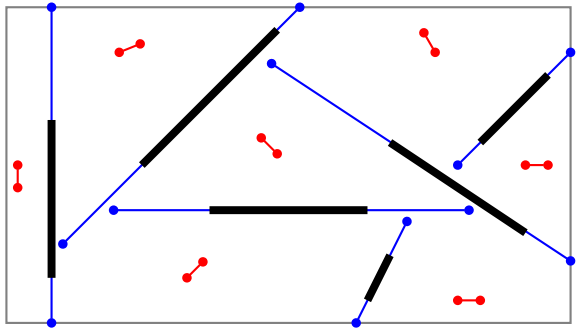}
\caption{The blue and red segments form a perfect matching with no disjoint compatible perfect matching.}
\label{GeneralOddExample}
\end{figure}

For a more general example, start with a perfect matching of $n$ black segments enclosed in a bounded convex region $C$. Extend each black segment, one by one, until it hits the boundary of $C$, or stop an $\epsilon$  distance from another segment, or extension of a segment. This gives a new perfect matching with $n$ blue segments. The blue segments form $n+1$ `regions' inside $C$. In the middle of each region insert a short red segment. The blue and red segments together form a perfect matching $M$ with $2n+1$ segments, as illustrated in Figure~\ref{GeneralOddExample}. No two red vertices are visible (for small enough $\epsilon$). So if $M$ has a disjoint compatible perfect matching, then every red vertex is paired with a blue vertex, which is impossible because there are $2n+2$ red vertices and $2n$ blue vertices. Thus $M$ has no disjoint compatible perfect matching. In fact, in the visibility graph of $V(M)$ minus $E(M)$, the red vertices form an independent set with more than half the vertices. Hence the visibility graph of $V(M)$ minus $E(M)$ has no (graph-theoretic) perfect matching, which in turn implies that $M$ has no disjoint compatible perfect matching

\section{Compatible Disjoint Matchings}
\label{sec:CompatibleDisjoint}

In this section we introduce a sequence of conjectures that imply the Compatible Matching Conjecture. Our general approach, given a perfect matching $M$, is to first compute an extension $L$ of $M$, then assign each vertex of $M$ to one of its two neighbouring cells in the convex subdivision formed by $L$, and then compute a perfect matching of the vertices assigned to each cell, the union of which is a perfect matching that is disjoint and compatible with $M$. The assignment of the vertices of $M$ to cells is modelled by an orientation of the edges of the dual multigraph.

\medskip\noindent\textbf{Extension Conjecture.} Every even perfect matching $M$ has an extension $L$, such that the associated dual multigraph $G$ admits an even orientation, with the property that whenever a vertex $v$ of $G$ has indegree $2$, the two incoming edges at $v$ do not arise from the same segment in $M$.
\smallskip

\begin{lemma}
\label{SS2DC}
The Extension Conjecture implies the Compatible Matching Conjecture.
\end{lemma}

\begin{proof}
Given a perfect matching $M$, let $L$ be an extension of $M$ with the properties stated in the Extension Conjecture. Let $G$ be the dual multigraph associated with $M$ and $L$. For each oriented edge $xy$ of $G$ corresponding to a vertex $v$ of $M$, \emph{assign} $v$ to the cell $y$. For each cell $y$, let $S_y$ be the set of vertices assigned to $y$. Since the orientation of $G$ is even, $|S_y|$ is even. Now $y$ is convex, and every vertex in $S_y$ is on the boundary of $y$. Thus $S_y$ is in convex position. Let $M_y$ be the subgraph of $M$ induced by $S_y$. Since no segment in $M$ intersects the interior of $y$, $M_y$ is a matching of $S_y$, and every segment in $M_y$ is on the convex hull of $S_y$. Moreover, by the final assumption in the Extension Conjecture, if $|S_y|=2$, then the two vertices are not adjacent in $M_y$. Thus Lemma~\ref{lem:Convex} is applicable to $S_y$ with the matching $M_y$. Thus $M_y$ has a disjoint compatible perfect matching $M'_y$. Let $M':=\cup_yM'_y$. Since each point is assigned to exactly one cell, $M'$ is a perfect matching. Since the cells are disjoint, and each $M'_y$ is non-crossing, $M'$ is non-crossing. Since the interior of every edge in $M'_y$ is contained in the interior of $y$, and every edge in $M$ only intersects $y$ at a vertex on the boundary, $M$ and $M'$ are compatible.
\end{proof}

\medskip\noindent\textbf{Two Subgraphs Conjecture.} Every even perfect matching $M$ has an extension $L$, such that the associated dual graph $G$ has an edge-partition into two subgraphs $G_1$ and $G_2$, such that each component of $G_1$ is even, each component of $G_2$ is even, and for every segment $vw$ of $M$, the edge of $G$ corresponding to $v$ is in a different subgraph from the edge of $G$ corresponding to $w$.
\smallskip

\begin{lemma}
\label{TS2TT}
The Two Subgraphs Conjecture implies the Extension Conjecture.
\end{lemma}

\begin{proof}
Given a perfect matching $M$, let $L$ be an extension of $M$ with the properties stated in the Two Subgraphs Conjecture. Since each component of $G_1$ and $G_2$ is even, by Lemma~\ref{lem:EvenOrientation}, each of $G_1$ and $G_2$ admit an even orientation. The union of the even orientations of $G_1$ and $G_2$ define an even orientation of $G$, such that if a vertex $x$ of $G$ has indegree $2$, then the two incoming edges at $x$ are both in $G_1$ or both in $G_2$, and thus arise from distinct segments in $M$. Hence the even orientation of $G$ satisfies the requirements of the Extension Conjecture. 
\end{proof}

\medskip\noindent\textbf{Two Trees Conjecture.} Every (even or odd) perfect matching $M$ has an extension $L$, such that the associated dual graph $G$ has an edge-partition into two trees, and for every segment $vw$ of $M$, the edge of $G$ corresponding to $v$ is in a different tree from the edge of $G$ corresponding to $w$.
\smallskip

\begin{lemma}
\label{lem:TT2TS}
The Two Trees Conjecture implies the Two Subgraphs Conjecture.
\end{lemma}

\begin{proof}
Let $M$ be a perfect matching with $n$ edges. Assuming the Two Trees Conjecture, $M$ has an extension $L$, such that the associated dual graph $G$ has an edge-partition into two trees $T_1$ and $T_2$. Now $G$ has $n+1$ vertices and $2n$ edges. Each tree has at most $n+1$ vertices, and thus has at most $n$ edges. Since $G$ has $2n$ edges, each tree has exactly $n$ edges. That is, each tree is a spanning tree of $G$. In the case that $M$ is even (which is assumed in the Two Subgraphs Conjecture), $n$ is even. Thus $T_1$ and $T_2$ are connected subgraphs, each with an even number of edges, as desired.
\end{proof}

Some notes on the Two Trees Conjecture are in order:

\begin{itemize}

\item \citet{Tutte61} and \citet{NW-JLMS61} independently characterised the multigraphs $G$ that contain two edge-disjoint spanning trees as those that have at least $2k-2$ cross-edges in every partition of $V(G)$ into $k$ parts. (A \emph{cross-edge} has endpoints in distinct parts of the partition.)\ 

\item It is easily seen that in every minimum counterexample to the Two Trees Conjecture, for every segment $r$, there exists segments $s$ and $t$, such that the line extending $r$ intersects $s$, and the line extending $t$ intersects $r$.

\item Motivated by the present paper, \citet{BDDHIST-CCCG07} made some progress on the Two Trees Conjecture. They proved that every perfect matching has an extension such that the associated dual multigraph $G$ is $2$-edge-connected, which is a necessary condition for $G$ to have the desired partition into two trees.

\end{itemize}

\section{Vertical-Horizontal Matchings}
\label{sec:VerticalHorizontal}

\begin{theorem}
\label{thm:VerticalHorizontal}
Every perfect matching $M$ consisting of vertical and horizontal segments satisfies the Two Trees Conjecture.
\end{theorem}

\begin{figure}[htb]
\centering\includegraphics{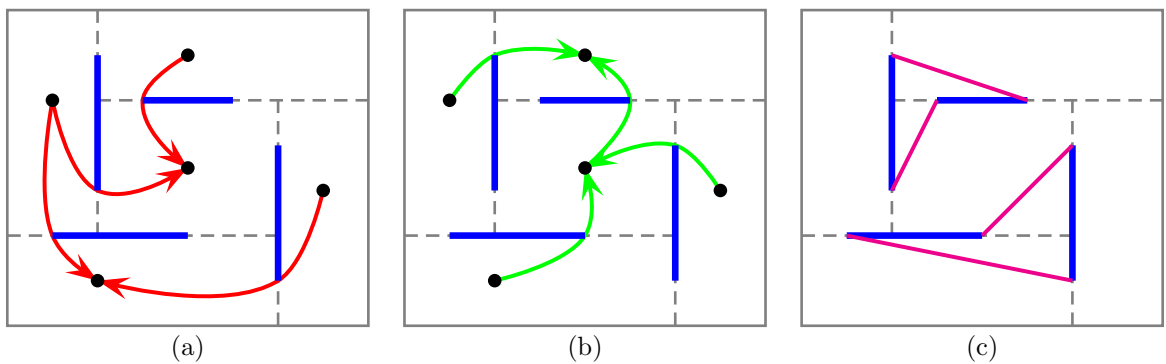}
\caption{For a perfect matching of horizontal and vertical segments: (a) the evenly oriented red spanning tree of the dual multigraph, (b) the evenly oriented green spanning tree, and (c) the compatible disjoint perfect matching determined by our algorithm.}
\label{fig:HoriVert} 
\end{figure}

\begin{proof}
Let $L$ be an extension of $M$ constructed as follows. First extend each horizontal segment in $M$ in both directions until it hits some vertical segment or goes to infinity. Then extend each vertical segment in $M$ in both directions until it hits some horizontal segment, an extension of some horizontal segment, or goes to infinity. Let $G$ be the dual multigraph associated with $M$ and $L$. 

Consider $G$ to be drawn in the plane with each vertex of $G$ in the interior of the corresponding cell of the convex subdivision formed by $L$. Each edge $xy$ of $G$ corresponding to a vertex $v$ of $M$ is drawn as a simple curve from $x$ through $v$ to $y$. Clearly $G$ can be drawn in this way without edge crossings. Moreover, if an edge of $G$ crosses $L$ then it crosses at a vertex of $M$. 

Colour the edges of $G$ \emph{red} and \emph{green} as follows. For each horizontal segment $vw$ in $M$, where $v$ is the left endpoint and $w$ is the right endpoint, colour the edge of $G$ through $v$ red, and colour the edge of $G$ through $w$ green. For each vertical segment $vw$ in $M$, where $v$ is the bottom endpoint and $w$ is the top endpoint, colour the edge of $G$ through $v$ red, and colour the edge of $G$ through $w$ green, as illustrated in Figure~\ref{fig:HoriVert}.


We claim that both the red and green subgraphs of $G$ are trees. Suppose on the contrary that $G$ has a monochromatic cycle $C$. Since the construction is symmetric between green and red, without loss of generality, $C$ is red. Now $C$ is a simple closed curve drawn without crossings in the plane.  Let \CC\ be the union of $C$ with its interior. 

First suppose that some horizontal segment of $M$ intersects \CC. Let $r$ be the topmost horizontal segment of $M$ that intersects \CC. ($r$ is uniquely determined since the vertices of $M$ are in general position.)\ Since the two edges of $G$ that pass through the endpoints of $r$ receive distinct colours, $C$ does not intersect both endpoints of $r$. If some endpoint of $r$ is in the interior of $C$ then let $v$ be that endpoint. Otherwise, since an edge of $G$ only intersects a segment of $M$ at the endpoint of that segment, some endpoint $v$ of $r$ is on $C$. In both cases, the extension of $r$ from $v$ enters the interior of $C$ and thus does not go to infinity because otherwise it would intersect $C$ at some point other than an endpoint of $r$. The extension of $r$ from $v$ is not blocked by some vertical extension because the horizontal segments were extended before the vertical segments. Thus the extension of $r$ from $v$ is blocked by some vertical segment $s$, and $s$ intersects \CC. Let $w$ be the top endpoint of $s$. Thus the Y-coordinate of $w$ is greater than that of $v$. Now $C$ does not pass through $w$ because the edge of $G$ through $w$ is green. The upward extension of $s$ does not go to infinity because otherwise it would intersect $C$ at some point other than an endpoint of $s$. Thus the upward extension of $s$ is blocked by some horizontal segment $t$, and $t$ intersects \CC. Thus $t$ is a horizontal segment of $M$ that intersects \CC\ and is higher than $r$. This contradiction proves that no horizontal segment of $M$ intersects \CC.
 
Every edge of $C$ passes through the endpoint of some segment $s$, in which case $s$ intersects \CC. Thus some segment of $M$ intersects \CC. Hence some vertical segment $r$ of $M$ intersects \CC. As in the previous case, if some endpoint of $r$ is in the interior of $C$ then let $v$ be that endpoint. Otherwise, some endpoint $v$ of $r$ is on $C$. In both cases, the extension of $r$ from $v$ does not go to infinity because otherwise it would intersect $C$ at some point other than an endpoint of $r$. Thus the extension of $r$ from $v$ is blocked by some horizontal segment $s$, and $s$ intersects \CC, which is a contradiction. 

Hence there is no monochromatic cycle in $G$. If $M$ has $n$ edges, then $G$ has $n+1$ vertices, the red subgraph has $n$ edges, and the green subgraph has $n$ edges. Every cycle-free graph with $n+1$ vertices and $n$ edges is a spanning tree. Thus the red subgraph is a tree and the green subgraph is a tree. By construction, for every segment $vw$ of $M$, the edge of $G$ passing through $v$ is in a different tree from the edge of $G$ passing through $w$. Thus the Two Trees Conjecture is satisfied.
\end{proof}

Theorem~\ref{thm:VerticalHorizontal} and Lemmas~\ref{SS2DC}--\ref{lem:TT2TS} imply:

\begin{corollary}
Every even perfect matching consisting of vertical and horizontal segments has a disjoint compatible perfect matching.
\end{corollary}

\section{Convex-Hull-Connected Matchings}
\label{CHC}

A set $M$ of pairwise disjoint segments is \emph{convex-hull-connected} if each segment has at least one endpoint on the boundary of the convex hull of $M$. This restriction has appeared in the context of augmenting a set of segments to form a simple polygon. Rappaport~et~al.~\citep{RIT-DCG90} gave an \Oh{n \log n} time algorithm to determine whether a set of convex-hull-connected segments admits an alternating polygon. Moreover, \citet{M-CGTA92} showed that every set of $n$ convex-hull-connected segments admits a circumscribing polygon, which can be constructed in \Oh{n \log n} time.

\begin{theorem}
\label{thm:CHC}
For every even convex-hull-connected perfect matching $M$, there is a perfect matching that is disjoint and compatible with $M$.
\end{theorem}

\begin{figure}[htb]
  \centering
  \includegraphics{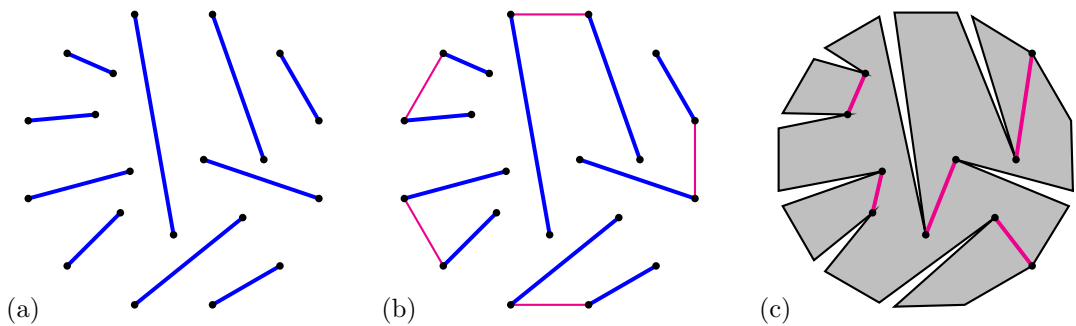}
  \caption{Proof of Theorem~\ref{thm:CHC}: (a) a convex-hull-connected set $M$ of segments, (b) the set $B$ of alternate gaps, (c) the polygon $P$ and matching $Q$.}
  \label{fig:CHC}
\end{figure}

\begin{proof}
We proceed by induction on the number of segments in $M$. A segment $vw$ in $M$ is a \emph{splitter} if $v$ and $w$ are non-consecutive points on the boundary of the convex hull (amongst the set of endpoints of segments in $M$).

First suppose that there is a splitter segment $vw$ in $M$. Of the sets of segments on the two sides of $vw$, one has an even non-zero number of segments, and the other has an odd number of segments. Group $vw$ with the odd-sized set. Thus $M$ is now partitioned into two even convex-hull-connected perfect matchings $M_1$ and $M_2$. By induction, there is a perfect matching $M_1'$ that is disjoint and compatible with $M_1$, and there is a perfect matching $M_2'$ that is disjoint and compatible with $M_2$. Hence $M_1'\cup M_2'$ is a perfect matching that is disjoint and compatible with $M$. 

Now assume $M$ has no splitter segment; refer to Figure~\ref{fig:CHC}. A \emph{gap} is an edge of the convex hull of $M$ that is not a segment in $M$. Since $M$ is even and there are no splitter segments, the number of gaps is even. Let $B$ be a set of alternate gaps on the convex hull. Thus $B$ forms a set of segments, such that for every segment $xy$ in $M$, exactly one of $x$ and $y$ is an endpoint of a segment in $B$. For each segment $xy$ with exactly one endpoint, say $x$, on the convex hull, let $W(xy)$ be an infinitesimally thick wedge centred at $y$ containing $xy$. Let $P$ be the polygon obtained from the convex hull of $M$ by removing each $W(xy)$. Thus every reflex vertex of $P$ is an endpoint of a segment in $M$ not intersecting $B$. Since $M$ is even and $B$ includes exactly one endpoint from each segment in $M$, the number of endpoints of segments in $M$ not intersecting $B$ is even. By Lemma~\ref{lem:reflex}, there is a perfect matching $Q$ of the set of endpoints of segments in $M$ not intersecting $B$, such that every segment in $Q$ is inside polygon $P$. Since every segment in $B$ is on the boundary of the convex hull, $B\cup Q$ is a perfect matching that is disjoint and compatible with $M$.  
\end{proof}

\section{The 4/5 Theorem}
\label{ThreeQuarters}

Given that the Compatible Matching Conjecture has remained elusive, it is natural to consider how large a disjoint compatible matching can be guaranteed.

\begin{theorem} 
\label{thm:FourFifths}
Let $S$ be a set of $2n$ points in the plane in general position, with $n$ even, and let $M$ be a perfect matching of $S$. Then there is a matching $M'$ of $S$ with at least $\fifth(4n-1)$ segments, such that $M$ and $M'$ are compatible and disjoint.
\end{theorem}

\begin{proof}
Without loss of generality, no segment is vertical.  
Fix a bounding box around the segments.
First extend each segment to the right (in any order).
Then extend each segment to the left (in any order).
We obtain a convex subdivision with $n+1$ faces.

Let $G$ be the corresponding dual multigraph. So $G$ has $n+1$ vertices. 
Colour each edge of $G$ that corresponds to a right endpoint \emph{red}.
Colour each edge of $G$ that corresponds to a left endpoint \emph{blue}.
Let $R$ and $B$ be the subgraphs of $G$, both with vertex set $V(G)$, 
respectively consisting of the red and blue edges. 
Each of $B$ and $R$ have $n+1$ vertices and $n$ edges.

We claim that $B$ is a spanning tree of $G$. To see this, consider the dual graph after the segments have been extended to the right, and before the segments have been extended to the left. At this point, the extensions do not form a cycle, and thus the dual graph has only one vertex. After extending the first segment to the left, the dual graph has two vertices, and one edge through a right endpoint. Thus the edges of the dual corresponding to the right endpoints form a spanning tree (a single edge) of the dual graph. As each subsequent extension of a segment to the left, the subgraph of the dual consisting of the edges through the right endpoints is obtained from the previous subgraph by a vertex splitting operation. It follows that after each left extension, the edges of the dual through the right endpoints form a spanning tree of the dual graph. Hence $B$ is a spanning tree of $G$. 

For each odd component $X$ of $R$, there is an edge $e$ in $X$, such that $X-e$ has no odd component. (\emph{Proof}: If $X$ has a leaf, then let $e$ be that edge. Otherwise $X$ has a cycle, and let $e$ be an arbitrary edge in the cycle.)\
Delete $e$ from $R$.
We are left with a subgraph $R'$ of $R$ with no odd component. 
Since $n$ is even, the one component of $B$ is even.
By construction, for every segment $vw$ of $M$, the edge of $G$ corresponding to $v$ is coloured differently from the edge of $G$ corresponding to $w$. 
Hence $B\cup R'$ satisfies the Two Subgraphs Conjecture.
By Lemma~\ref{TS2TT}, there is a partial matching $M'$ of $S$ that is compatible and disjoint with $M$, and the number of segments in $M'$ equals half the number of edges in $B\cup R'$, which is $2n$ minus the number of odd components in $R$.
Lemma~\ref{lem:OddComps} below bounds the number of odd components in a planar graph. This result applied to $R$ (which has $n+1$ vertices, $n$ edges, and thus has some  component not isomorphic to $K_2$) implies that $R$ has at most $\fifth(3(n+1)-n-1)=\frac25(n+1)$ odd components. 
Hence $M'$ has at least $\half(2n-\frac25(n+1))=\tfrac{1}{5}(4n-1)$ segments.
\end{proof}

\begin{lemma}
\label{lem:OddComps}
Let $f(G)$ be the number of odd components in a graph $G$. Then every planar graph $G$ with $n$ vertices and $m$ edges has $f(G)\leq\fifth(3n-m)$, with equality only if every component of $G$ is $K_2$.
\end{lemma}

\begin{proof}
We proceed by induction on the number of components in $G$.

For the base case, suppose that $G$ has one component.
If $n=1$ then $f(G)=0<\frac{3}{5}=\fifth(3n-m)$.
If $n=2$ then $f(G)=1=\fifth(3n-m)$.
If $n\geq3$ then $f(G)\leq 1<\frac{6}{5}\leq \fifth(3n-(3n-6))\leq \fifth(3n-m)$.
Now assume that $G$ has at least two components. 

Suppose that $G$ has an isolated vertex $v$.
By induction, $f(G)=f(G-v)\leq\fifth(3(n-1)-m)<\fifth(3n-m)$.
Now assume that $G$ has no isolated vertices.

Suppose that $G$ has a component $H=K_2$ . Then $f(G)=f(G-H)+1\leq\fifth(3(n-2)-(m-1))+1=\fifth(3n-m)$. Moreover, suppose that $f(G)=\fifth(3n-m)$. Then $f(G-H)=\fifth(3(n-2)-(m-1))$, and by induction, every component of $G-H$ is $K_2$, which implies that every component of $G$ is $K_2$.
Now assume that $G$ has no $K_2$ component.

Let $X$ be a component of $G$ with $p$ vertices and $q$ edges.
By induction, $f(G)\leq 1+f(G-X)\leq 1+\fifth(3(n-p)-(m-q))
=\fifth(3n-m+q-3p+5)<\fifth(3n-m)$ since $q\leq 3n-6$.
\end{proof}

\begin{figure}[htb]
\begin{center}
\includegraphics[scale=0.25]{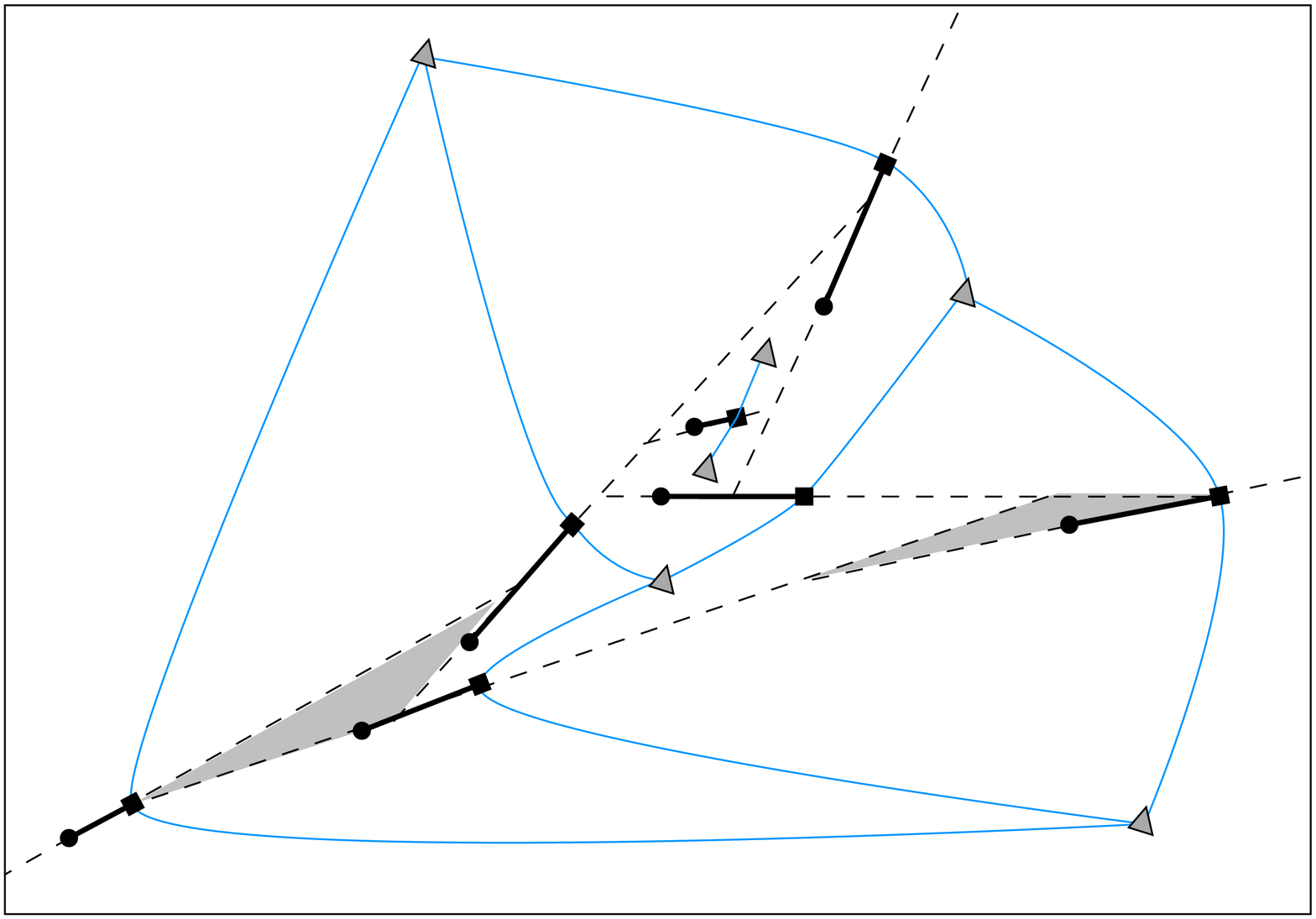}
\includegraphics[scale=0.325]{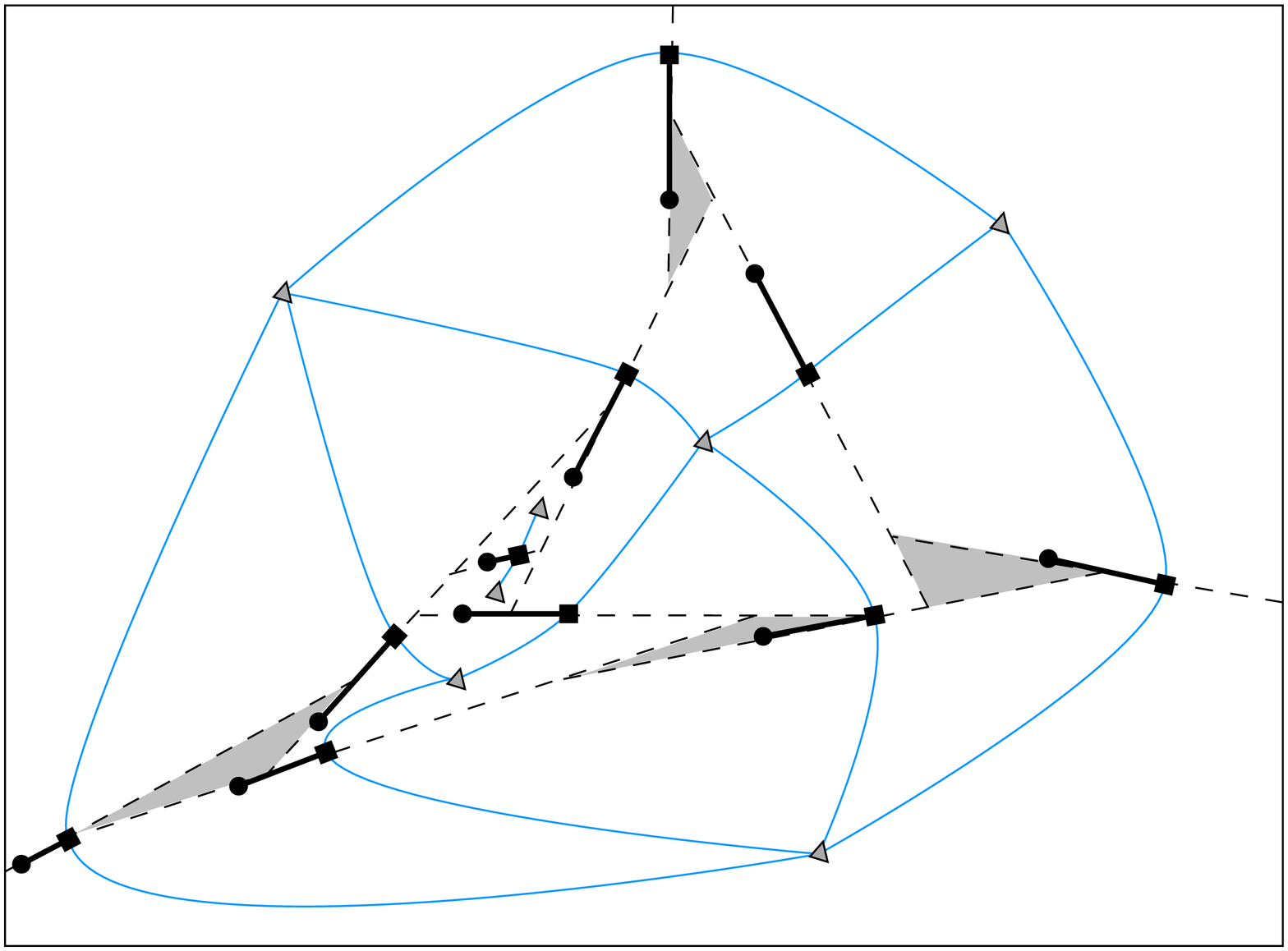}
\caption{Consider the seven segments in the figure on the left. It is possible to extend the segments (all right endpoints first), so that the graph $R$ in Theorem~\ref{thm:FourFifths} (the spanning subgraph of the dual graph consisting of the edges passing through the right endpoints) has two components, $K_4$ and $K_2$. Placing an additional segment in each of the two shaded regions creates two additional $K_2$ components in $R$. Then $R$ has one component that is maximal planar, and every other component is $K_2$. By repeatedly adding two more segments as shown in the figure on the right, we obtain arbitrarily large sets of segments such that the graph $R$ has 
one maximal planar component, and every other component is $K_2$.} 
\label{fig:construct}
\end{center}
\end{figure}

We now show that the analysis of the algorithm in the proof of Theorem~\ref{thm:FourFifths} is tight. First note that if an $n$-vertex $m$-edge planar graph $G$ has one component that is maximal planar on an odd number of vertices, and every other component is $K_2$, then $f(G)=\fifth(3n-m-1)$. Figure~\ref{fig:construct} shows a set of segments such that by applying the algorithm in the proof of Theorem~\ref{thm:FourFifths}, the obtained graph $R$ has one component that is maximal planar on an odd number of vertices, and every other component is $K_2$. It follows that for this set of segments, the algorithm in Theorem~\ref{thm:FourFifths} will produce a matching with $\fifth(4n-1)$ segments.

\section{Matchings with Crossings}
\label{WithCrossings}

Now we relax the Compatible Matching Conjecture by allowing crossings.

\begin{theorem}
\label{thm:ConCrossings}
Let $M$ be an even perfect matching with no vertical segment. Let $L$ be the set of left endpoints of $M$, and let $R$ be the set of right endpoints of $M$. Then there is a perfect matching $M_L$ of $L$, and a perfect matching $M_R$ of $R$, such that no edge in $M$ crosses an edge in $M_L\cup M_R$ (but an edge in $M_L$ might cross an edge in $M_R$).
\end{theorem}

\begin{proof}
Let $C$ be a convex polygon bounding $M$. In a similar fashion to the second proof of Lemma~\ref{lem:HalfPlane} and the proof of Theorem~\ref{thm:CHC}, extend each segment of $M$ by an infinitesimally thickened ray from its left endpoint. Removing the thickened rays from the interior of $C$, we obtain a polygon whose reflex vertices are the right endpoints of the segments in $M$. Since $M$ is even, by Lemma~\ref{lem:reflex} with $R=S$, there is a perfect matching $M_R$ of $R$ such that $M_R\cup M$ is non-crossing. The perfect matching $M_L$ is obtained similarly.
\end{proof}

\begin{corollary}
\label{cor:ConCrossings}
Let $M$ be an even perfect matching. Let $G$ be the visibility graph of $V(M)$ minus $E(M)$. Then $G$ contains a graph-theoretic perfect matching (which possibly has crossings, but is the union of two non-crossing matchings).
\end{corollary}

Note that the assumption that $M$ is even is needed in Corollary~\ref{cor:ConCrossings}---because of the instance in Figure~\ref{OddExample} for example.


\def\soft#1{\leavevmode\setbox0=\hbox{h}\dimen7=\ht0\advance \dimen7
  by-1ex\relax\if t#1\relax\rlap{\raise.6\dimen7
  \hbox{\kern.3ex\char'47}}#1\relax\else\if T#1\relax
  \rlap{\raise.5\dimen7\hbox{\kern1.3ex\char'47}}#1\relax \else\if
  d#1\relax\rlap{\raise.5\dimen7\hbox{\kern.9ex \char'47}}#1\relax\else\if
  D#1\relax\rlap{\raise.5\dimen7 \hbox{\kern1.4ex\char'47}}#1\relax\else\if
  l#1\relax \rlap{\raise.5\dimen7\hbox{\kern.4ex\char'47}}#1\relax \else\if
  L#1\relax\rlap{\raise.5\dimen7\hbox{\kern.7ex
  \char'47}}#1\relax\else\message{accent \string\soft \space #1 not
  defined!}#1\relax\fi\fi\fi\fi\fi\fi} \def\Dbar{\leavevmode\lower.6ex\hbox to
  0pt{\hskip-.23ex \accent"16\hss}D}

\end{document}